\documentclass[12pt]{article}

\usepackage{amssymb,amsthm,amsfonts,amsmath,amscd,verbatim}
\usepackage{xy}
\xyoption{all}
\title{On closed rational functions in several variables}
\author{A.~P.~Petravchuk\footnote{Kiev Taras Shevchenko University, Faculty of
  Mechanics and Mathematics, 64, Volodymyrska street, 01033 Kyiv,
  Ukraine, e-mail: aptr@univ.kiev.ua}, O.~G.~Iena\footnote{ Kiev Taras Shevchenko
   University and
 Technische
    Universit\"at Kaiserslautern, Fachbereich Mathematik,
Postfach 3049, 67653 Kaiserslautern, Germany, e-mail:
yena@mathematik.uni-kl.de}}
\date{}

\newtheorem{tr}{Theorem}
\newtheorem*{tr*}{Theorem}

\newtheorem{lemma}{Lemma}

\newtheorem*{pr*}{Proposition}

\newtheorem{cor}{Corollary}
\theoremstyle{definition}
\newtheorem{df}{Definition}

\newtheorem*{df*}{Definition}
\newtheorem*{not*}{Notation}
\newtheorem*{ex}{Example}

\newtheorem{rem}{Remark}

\newtheorem*{rem*}{Remark}

\DeclareMathOperator\im{Im}

\DeclareMathOperator\trdeg{tr.deg}
\def\differential{d}
\renewcommand\Im\im
\renewcommand\d\differential

\def\bb#1{\mathbb #1}

\def\dd#1#2{\frac{\partial #1}{\partial#2} }
\newcommand{\rpoly}[1]{\k(x_1,\dots, x_{#1})}
\newcommand{\poly}[1]{\k[x_1,\dots, x_{#1}]}

\def\eqdef{=:}

\def\k{\mathbb K}

\let\le\leq
\let\ge\geq
\let\star *
\let\subset\subseteq

\begin{document}
\maketitle
\begin{abstract}
Let $\k=\bar \k$ be a field of characteristic zero. An element
$\varphi\in \rpoly{n}$ is called a closed rational function if the
subfield $\k(\varphi)$ is algebraically closed in the field
$\rpoly{n}$. We prove that a rational function $\varphi=f/g$ is
closed if $f$ and $g$ are algebraically independent and at least
one of them is irreducible. We also show that the rational
function $\varphi=f/g$
is closed if and only if the pencil $\alpha f+\beta g$ contains
only finitely many reducible hypersurfaces. Some sufficient
conditions for a polynomial to be irreducible are given.
\end{abstract}
\section{Introduction}
Closed polynomials, i. e., polynomials $f\in \poly{n}$ such that
the subalgebra $\k[f]$ is integrally closed in $\poly{n}$, were
studied by many authors (see, for example, \cite{Now},
\cite{Stein}, \cite{Najib},\cite{Arzh} ). A rational analogue of a
closed polynomial is a rational function $\varphi$ such that the
subfield $\k(\varphi)$ is algebraically closed in the field
$\rpoly{n}$, such a rational function will be called a closed one.
Although there are algorithms to determine whether a given
rational function is closed (see~\cite{Oll}) it is interesting to
study closed rational functions  more detailed.

We give the following sufficient condition for a rational function
to be closed. Let $\varphi=f/g\in \rpoly{n}$, $f$ and $g$ are coprime,
algebraically independent and at least one of polynomials $f$ and $g$
is irreducible. Then $\varphi$ is a closed rational function
(Theorem~\ref{th:condition_for_closed}).

Using some results of of J.M. Ollagnier~\cite{Oll} about Darboux
polynomials we prove that a rational function
$\varphi=f/g\in\rpoly{n}\setminus \k$ is closed if and only if the
pencil $\alpha f+\beta g$ of hypersurfaces contains only finitely
many reducible hypersurfaces (Theorem~\ref{th:criterion}). We also
study
products  of irreducible polynomials.

 Notations in the paper are standard.
For a rational function $F(t)\in \k(t)$  of the form
$F(t)=\frac{P(t)}{Q(t)}$ with coprime polynomials $P$ and $Q$ the
degree is $\deg F=\max(\deg P, \deg Q)$. The ground field $\k$ is
algebraically closed  of characteristic $0$.

\section{Closed rational functions in several variables}

\begin{lemma}\label{lemma:alg_dep}
For a rational functions $\varphi, \psi\in \rpoly{n}\setminus \k$ the
following conditions are equivalent.

1) $\varphi$ and $\psi$ are algebraically dependent over $\k$;

2) the rank of Jacobi matrix
$
J(\varphi, \psi)=
\begin{pmatrix}
\dd{\varphi}{x_1}&\dots&\dd{\varphi}{x_n}\\
\dd{\psi}{x_1}&\dots&\dd{\psi}{x_n}
\end{pmatrix}
$ is equal to $1$;

3) for differentials $\d{\varphi}$ and $\d{\psi}$ of functions $\varphi$ and $\psi$
respectively it holds $\d{\varphi}\wedge\d{\psi}=0$;

4) there exists $h\in \rpoly{n}$ such that $\varphi=F(h)$ and $\psi=G(h)$ for
 some $F(t), G(t)\in \k(t)$.
\end{lemma}
\begin{proof}
The equivalence of 1) and 2) follows from~\cite{Hodge}, Ch.III, \S 7, Th.
III. The equivalence of 2) and 3) is obvious. Since 2) clearly
follows from 4), it remains to show that 1) implies 4). Let
$\varphi$ and $\psi$ be algebraically dependent. Then obviously
$\trdeg_{\k}\k(\varphi, \psi)=1$. By Theorem of Gordan (see for
example~\cite{Schinzel}, p.15) $\k(\varphi, \psi)=\k(h)$ for some
rational function $h$ and therefore $\varphi=F(h)$ and $\psi=G(h)$
for  some $F(t), G(t)\in \k(t)$.
\end{proof}

\begin{df}
We call a rational function $\varphi \in \rpoly{n}\setminus \k$  closed if the
subfield $\k(\varphi)$ is algebraically closed in $\rpoly{n}$.

A rational function $\tilde \psi$ is called generative for a rational function $\psi$
if $\tilde \psi$ is closed and $\psi\in \k(\tilde \psi)$.
\end{df}

\begin{lemma}\label{lemma:generative}
1) For a rational function $\varphi\in \rpoly{n}\setminus \k$ the
following conditions are equivalent.

a) $\varphi$ is closed;

b) $\k(\varphi)$ is a maximal element in the partially ordered (by
inclusion) set of subfields of $\rpoly{n}$ of the form $\k(\psi)$,
$\psi\in \rpoly{n}\setminus \k$.

c) $\varphi$ is non-composite rational function, i. e., from the
equality $\varphi=F(\psi)$, for some rational functions $\psi\in
\rpoly{n}\setminus \k$ and $F(t)\in \k(t)$, it follows that $\deg
F=1$.

2) For every rational function $\varphi\in \rpoly{n}\setminus\k$
there exists a generative rational function $\tilde \varphi$. If
$\tilde \varphi_1$ and $\tilde \varphi_2$ are two generative
rational functions for $\varphi$, then $\tilde
\varphi_2=\frac{a\tilde\varphi_1+b}{c\tilde\varphi_1+d}$ for some
$a, b, c, d\in \k$ such that $ad-bc\not= 0.$

\end{lemma}
\begin{proof}
1) a)$\Rightarrow$ b). Suppose that the rational function $\varphi$ is closed
and $\k(\varphi)\subset\k(\psi)$ for some $\psi\in \rpoly{n}$. The element $\psi$
is algebraic over over $\k(\varphi)$ and therefore by the definition of closed rational
functions we have $\psi\in \k(\varphi)$. Thus  $\k(\varphi)$ is a maximal element in the
set of all one-generated subfields of $\rpoly{n}$.

b)$\Rightarrow$ a). If $\k(\varphi)$ is a maximal one-generated subfield of $\rpoly{n}$,
then $\k(\varphi)$ is algebraically closed in  $\rpoly{n}$. Indeed, if $f$ is algebraic
over $\k(\varphi)$, then $\trdeg_\k(\varphi, f)=1$ and by Theorem
of Gordan $\k(\varphi, f)=\k(\psi)$ for some rational function $\psi$.
But then $\k(\psi)=\k(\varphi)$ and $f\in \k(\varphi)$.

The equivalence of b) and c) is obvious.

2) The subfield $\k(\varphi)$ is contained in some maximal one-generated
subfield $\k(\tilde\varphi)$, which is algebraically closed in $\rpoly{n}$ by
the part 1) of this Lemma. Therefore $\tilde\varphi$ is a generative rational
function for $\varphi$.

Let  $\tilde\varphi_1$ and $\tilde\varphi_2$ be two generative
rational functions for $\varphi$. Then $\tilde\varphi_1$ and
$\tilde\varphi_2$ are algebraic over the field $\k (\varphi )$ and
therefore $\trdeg _{\k}\k (\varphi ,\tilde\varphi_1,
\tilde\varphi_2)=1$. In particular, then the  rational functions
$\tilde\varphi_1$ and $\tilde\varphi_2$ are algebraically
dependent.

  By Lemma~\ref{lemma:alg_dep} one obtains $\tilde\varphi_1\in
\k(\psi)$, $\tilde\varphi_2\in \k(\psi)$ for some rational
function $\psi$. Since both $\tilde\varphi_1$ and
$\tilde\varphi_2$ are closed, we get
$\k(\tilde\varphi_1)=\k(\psi)=\k(\tilde\varphi_2)$. But there
exists a fractional rational transformation  $\theta$ of the field
$\k(\psi)$ such that $\theta(\tilde\varphi_2)=\tilde\varphi_1$.
Therefore,
$\tilde\varphi_2=\frac{a\tilde\varphi_1+b}{c\tilde\varphi_1+d}$,
for some $a, b, c, d\in \k ,$  $ad-bc \neq 0$.
\end{proof}

\begin{rem}
Note that algebraically dependent rational functions have the same set of generative
functions. This follows from Lemma~\ref{lemma:alg_dep} and Lemma~\ref{lemma:generative}.
\end{rem}

\begin{rem}
\label{rem:closed_poly_are_closed_rat} Let $f\in \poly{n}\setminus \k .$ By Lemma 3 from~\cite{Arzh},
 the subfield $\k (f)$ is algebraically closed if and only
if the polynomial $f$ is closed. So, the polynomial $f$ is closed if
and only if $f$ is closed as a rational function.
\end{rem}

\begin{tr}\label{th:condition_for_closed}
Let polynomials $f, g\in \poly{n}$ be coprime and algebraically
independent. If at least one of them is  irreducible, then the
rational function $\varphi=\frac{f}{g}$ is closed.
\end{tr}
\begin{proof}
Without loss of generality we can assume  that $f$ is irreducible.
By Lemma~\ref{lemma:generative} there exists a generative rational
 function $\psi=\frac{p}{q}$ for $\varphi$, where $p$ and $q$ are coprime
 polynomials. Then $\varphi=\frac{P(\psi)}{Q(\psi)}$ for some coprime polynomials $P(t),
 Q(t)\in \k[t]$.

Let $P(t)=a_0(t-\lambda_1)\dots(t-\lambda_m)$ and
$Q(t)=b_0(t-\mu_1)\dots(t-\mu_l)$, $\lambda_i, \mu_j\in \k$ be
 the decompositions of $P(t)$ and $Q(t)$ into irreducible factors. Then
\[
\varphi=\frac{f}{g}=\frac{a_0(\frac{p}{q}-\lambda_1)\dots(\frac{p}{q}-\lambda_m)}{b_0(\frac{p}{q}-\mu_1)\dots(\frac{p}{q}-\mu_l)}=
\frac{ a_0(p-\lambda_1 q)\dots(p-\lambda_m q)q^{l-m}}{b_0(p-\mu_1 q)\dots(p-\mu_l q)}.
\]
and we obtain
\[
b_0f(p-\mu_1 q)\dots(p-\mu_l q)=a_0g(p-\lambda_1 q)\dots(p-\lambda_m q)q^{l-m}.\leqno{(\star)}
\]

Note, as $\lambda_i\neq \mu_j$, the polynomials $p-\lambda_iq$
and $p-\mu_jq$ are coprime for  all $i=\overline{1, l}$ and $j=\overline{1,  m}$.
 Moreover, since $p$ and $q$ are coprime, it is clear that $q$ is coprime with
 polynomials of the form $p+\alpha q$, $\alpha\in \k$.

Note also that  $p-\beta q\not \in\k$. Indeed, if $p-\beta
q=\xi\in \k$ for some $\beta, \xi\in \k$, then $p=\xi+\beta q$ and
\[
\varphi=\frac{f}{g}=
\frac{ a_0(\xi+(\beta-\lambda_1 )q)\dots(\xi+(\beta-\lambda_m )q)q^{l-m}}{b_0(\xi+(\beta-\mu_1 q)\dots(\xi+(\beta-\mu_l) q)},
\]
which means that $f$ and $g$ are algebraically dependent, which contradicts our  assumptions.

So from $(\star)$ we conclude that $f$ is divisible by all polynomials
 $(p-\lambda_iq)$. Since $f$ is irreducible, taking into account the above
 considerations  we conclude that $m=1$ and $f=a(p-\lambda_1 q)$, $a\in \k^\star$.
Therefore, from $(\star)$ we obtain
\[
b_0a(p-\lambda_1q)(p-\mu_1 q)\dots(p-\mu_l q)=a_0g(p-\lambda_1 q)q^{l-1}
\] and after reduction
\[
b_0a(p-\mu_1 q)\dots(p-\mu_l q)=a_0gq^{l-1}.
\]
Since $q$ is coprime with $(p-\mu_jq)$, we get $l=1$ and finally
$a_0g=b_0a(p-\mu_1q)$. We obtain
\[
\varphi=\frac{f}{g}=\frac{a_0 (p-\lambda_1
q)}{b_0(p-\mu_1q)}=\frac{a_0(\frac{p}{q}-\lambda_1)}{b_0(\frac{p}{q}-\mu_1)}=\frac{a_0(\psi-\lambda_1)}{b_0(\psi-\mu_1)}.
\]
One concludes that $\k(\varphi)=\k(\psi)$, which means that $\varphi$ is a
closed rational function.
\end{proof}

\section{Rational functions and pencils of hypersurfaces}

In this section we give a characterization of closed rational
functions. While proving Theorem ~\ref{th:criterion} we use an
approach from the paper of  J.M.Ollagnier~\cite{Oll} connected
with Darboux polynomials. Recall some notions and terminology (see
also \cite{Now1}, pp.22-24). If $\delta$ is a derivation of the
polynomial ring $\poly{n}$, then a polynomial $f$ is called a
Darboux polynomial for $\delta$ if $\delta (f)=\lambda f$ for some
polynomial $\lambda$ (not necessarily $\lambda \in \k$). The
polynomial $\lambda $ is called the cofactor for $\delta$
corresponding to the Darboux polynomial $f$ (so, $f$ is a
polynomial eigenfunction for $\delta$ and $\lambda$ is the
corresponding eigenvalue).

Further, for a rational function $\varphi =\frac{f}{g}\in \k
(x_{1},\ldots , x_{n})\setminus \k$ one can define a (vector)
derivation $\delta _{\varphi}=gdf-fdg : \poly{n}\to \Lambda
^{2}\poly{n}$ by the rule $\delta _{\varphi}(h)=dh\wedge
(gdf-fdg).$ For such a derivation $\delta _{\varphi}$ a polynomial
$h$ is called a Darboux polynomial  if all coefficients of the
$2-$form $dh\wedge (gdf-fdg)$ are divisible by the polynomial $h$,
i. e., $dh\wedge (gdf-fdg)=h\cdot \lambda$ for some $2-$form
$\lambda$, which is called a cofactor for the derivation $\delta
_{\varphi}.$ Note that every divisor of the Darboux polynomial $h$
is also a Darboux polynomial for $\delta _{\varphi}$ (see, for
example, \cite{Now1}, p.23). It is easy to see that the polynomial
$\alpha f+\beta g$ is a Darboux polynomial for the derivation
$\delta _{\varphi}$ and therefore every divisor of the polynomial
$\alpha f+\beta g$ is a Darboux polynomial of the derivation
$\delta _{\varphi}=gdf-fdg.$

\begin{tr}\label{th:criterion}
Let  polynomials $f,g\in \poly{n}$ be coprime and let at least one
of them be a  non-constant polynomial. Then the rational  function $\varphi =\frac{f}{g}$ is
closed  if and only if all but finitely many hypersurfaces in the
pencil $\alpha f+\beta g$ are irreducible.
\end{tr}
\begin{proof}
Let $\varphi=\frac{f}{g}$ be closed. Suppose that the pencil
$\alpha f+\beta g$ contains infinitely many reducible
hypersurfaces. Let $\{\alpha_i f+\beta_i g\}_{i\in \bb N}$,
$(\alpha_i:\beta_i)\neq (\alpha_j:\beta_j)$ for $i\neq j$ as
points of $\bb P^{1}$, be an infinite sequence of (different)
reducible hypersurfaces. For each $i$ take one irreducible factor
$h_i$ of $\alpha_i f+\beta_i g$.

By the above remark, all polynomials $h_i$ are Darboux polynomials
for $\delta_\varphi$ and $\deg h_i<\deg \varphi\eqdef k$. By
Corollary 5 from~\cite{Oll} there exist finitely many cofactors of
$\delta_\varphi$ that correspond to Darboux polynomials $h_i$
(degrees of $h_i$ are bounded). Therefore, there exist polynomials
$h_i$ and $h_j$ such that $\delta_\varphi(h_i)=\lambda h_i$ and
$\delta_\varphi(h_j)=\lambda h_j$ for some cofactor $\lambda\in
\bigwedge^2\poly{n}$. This implies
$\delta_\varphi(\frac{h_i}{h_j})=0$ and thus $\d
(\frac{f}{g})\wedge\d(\frac{h_i}{h_j})=\frac{1}{g^{2}}\delta_\varphi(\frac{h_i}{h_j})=0$
(see \cite{Oll}).  By Lemma~\ref{lemma:alg_dep}, the rational
functions $\varphi=\frac{f}{g}$ and $\frac{h_i}{h_j}$ are
algebraically dependent. As $\varphi$ is closed,
$\frac{h_i}{h_j}=F(\varphi)$ for some $F(t)\in \k(t)$ and
therefore $\deg\frac{h_i}{h_j}=\deg F\deg\varphi$ (see for
example~\cite{Oll}). But this is impossible since $\deg
\frac{h_i}{h_j}<\deg\varphi$ . Therefore, all but finitely many
hypersurfaces in $\alpha f+\beta g$ are irreducible.

Let now $\alpha_0 f+\beta _0 g$ be an irreducible hypersurface
from the pencil $\alpha f+\beta  g$. Consider the case when  $f$
and $g$ are algebraically independent. One can assume without loss
of generality $\alpha_0\neq 0$. Then
 $\alpha_0 f+\beta _0 g$ and $g$ are algebraically independent
 as well. (If $\alpha_0=0$, then $\beta_0\neq 0$ and
  polynomials $f$ and $\alpha_0f +\beta_0 g$ are algebraically independent).
  Therefore, since $\alpha_0 f+\beta_0 g$ and $g$  are coprime,
  by Theorem~\ref{th:condition_for_closed} the rational
  function $\psi=\frac{\alpha_0 f+\beta _0 g}{g}$ is closed.
   Therefore, $\varphi=\frac{f}{g}=\alpha_0^{-1}(\psi-\beta _0)$,
    which proves that $\varphi$ is a closed rational function.

Let now $f$ and $g$ be algebraically dependent.
Then $f=F(h)$ and $g=G(h)$ for a common generative polynomial
 function $h$ and polynomials $F(t), G(t)\in \k[t]$
 (see Remark~\ref{rem:closed_poly_are_closed_rat}).
 Let $(1:\beta_1)\neq (1:\beta_2)$ be two different
 points in $\bb P^1$ such that $f+\beta_i g$ is
 irreducible for $i\in\{1,2\}$. Then $f+\beta_i g=F(h)+\beta_i G(h)$
 is irreducible for  $i\in\{1,2\}$. In particular this
 means that $\deg(F(t)+\beta_i G(t))=1$, i. e.,
$F(t)+\beta_i G(t)=a_it+b_i$, $a_i, b_i\in \k$, $a_i\neq 0$. Then
since $\beta_1\neq \beta_2$, we conclude that $F(t)=at+b$ and
$G(t)=ct+d$ for some $a,b,c,d\in \k$. So
$\varphi=\frac{f}{g}=\frac{ah+b}{ch+d}$. As at least one of $f$
and $g$ is non-constant and since $f$ and $g$ are coprime, we
conclude that  $\k(\varphi)=\k(h)$. Therefore, since $\k(h)$ is an algebraically closed subfield of the field $\rpoly{n}$, $\varphi=\frac{f}{g}$ is a closed rational function.
\end{proof}

\begin{rem}
Note, in order to show that $\varphi=\frac{f}{g}$ is closed in
Theorem~\ref{th:criterion} it is enough to have two different
irreducible hypersurfaces in the pencil $\alpha f+\beta g$.
One irreducible hypersurface $\alpha_0 f+\beta_0 g$ is enough
provided $f$ and $g$ are algebraically independent.
\end{rem}

\begin{rem}
We also reproved a weak version (we do not give an estimation) of
the next result of W. Ruppert (see~\cite{Ruppert}, Satz 6).

If $f$ and $g$ are algebraically independent polynomials and the
pencil $\alpha f+\beta g$ contains at least one irreducible hypersurface,
then all but finitely many hypersurfaces in $\alpha f+\beta g$ are irreducible.
\end{rem}

\begin{rem}\label{clos}
If a polynomial  $f\in \poly{n}$ is non-constant then by
Theorem~\ref{th:criterion} and
Remark~\ref{rem:closed_poly_are_closed_rat} $f$ is closed if and
only if for all but finitely many $\lambda \in \k$ the polynomial
$f+\lambda$ is irreducible. This result is well-known (it can be
proved by using  the first Bertini theorem), see, for example,
\cite{Schinzel}, Corollary 3.3.1.
\end{rem}

Using   Remark~\ref{clos} one can prove that any polynomial $f\in
\poly{n}\setminus \k$ can be written in the form $f=F(h)$ for some
polynomial $F(t)\in \k[t]$ and irreducible polynomial $h$. Similar
statement holds for rational functions.

\begin{cor}\label{cor:irreducible}
A rational function $\frac{f}{g}\in \rpoly{n}\setminus \k$ can be
written in the form $\frac{f}{g}=F(\varphi)$, $F(t)\in \k(t)$, for
some rational function $\varphi=\frac{p}{q}$ such that polynomials
$p$ and $q$ are irreducible.
\end{cor}
\begin{proof}
Let $\frac{p_1}{q_1}$ be a generative function for $\frac{f}{g}$.
As $\frac{p_1}{q_1}$ is closed, by Theorem~\ref{th:criterion}
the pencil $\alpha p_1+\beta q_1$ contains two different
irreducible hypersurfaces $p=\alpha_1 p_1+\beta_1 q_1$ and
$q=\alpha_2 p_1+\beta_2 q_1$, i. e., $(\alpha_1:\beta_1)\neq (\alpha_2:\beta_2)$.
Since the pencils $\alpha p+\beta q$ and $\alpha p_1+\beta q_1$ are equal,
 and since in the pencil $\alpha p_1+\beta q_1$ all but finitely many
 hypersurfaces are irreducible, we conclude that $\frac{p}{q}$ is closed
 and is a generative function for $\frac{f}{g}$.
\end{proof}
\begin{rem}
Under conditions of Corollary~\ref{cor:irreducible} polynomials $p$ and $q$ can be chosen of the same degree.
\end{rem}
\begin{cor}
Let $\k\subset L\subset \rpoly{n}$ be  an algebraically closed subfield in $\rpoly{n}$. Then it is possible to choose generators of $L$ in the form $\frac{p_1}{q_1},\dots, \frac{p_m}{q_m}$, where $p_i$ and  $q_i$ are irreducible polynomials.
\end{cor}

\begin{tr}\label{tr:reducible}
Let polynomials $f, g\in \poly{n}$ be coprime and algebraically
independent. Then the rational function $f/g$ is not closed if and
only if there exist  algebraically independent  irreducible
polynomials $p$ and $q$
 and a positive integer $k\ge 2$ such that $f=(\alpha_1
p+\beta_1 q)\dots(\alpha_k p+\beta_k q)$ and $g=(\gamma_1
p+\delta_1 q)\dots(\gamma_k p+\delta_k q)$ for some
$(\alpha_i:\beta_i), (\gamma_j:\delta_j)\in \bb P_1$, with
$(\alpha_i:\beta_i)\neq(\gamma_j:\delta_j)$, $i,j=\overline{1,k}$.
\end{tr}
\begin{proof}
Suppose $\frac{f}{g}$ is not closed. Take its generative function
$\frac{p}{q}$ with irreducible polynomials $p$ and $q$. This is
 possible by Corollary~\ref{cor:irreducible}. Then $\frac{f}{g}=F(\frac{p}{q})$
 for some rational function $F(t)\in \k(t)$ with $\deg F(t)=k\ge 2$. Note that the polynomials
 $p$ and $q$ are algebraically independent because in other case the  polynomials $f$ and $g$
 were algebraically dependent which contradicts to our assumptions.
 Write
\[
F(t)=\frac{a_0(t-\lambda_1)\dots(t-\lambda_s)}{b_0(t-\mu_1)\dots(t-\mu_r)}
\] with $\lambda_i\neq \mu_j$, i. e., with coprime nominator
and denominator. It is clear that $k=\deg F(t)=\max\{s, r\}$.
After substitution of $t$ by $\frac{p}{q}$ we obtain
\[
\frac{f}{g}=\frac{a_0(p-\lambda_1 q)\dots(p-\lambda_s q)q^{r-s}}{b_0(p-\mu_1q)\dots(p-\mu_rq)}.
\]
Put $(\alpha_i:\beta_i)=(1:-\lambda_i)$ for $i=\overline{1,s}$, $(\gamma_j:\delta_j)=(1:-\mu_j)$ for $j=\overline{1, r}$. If $r\le s$, then put $(\gamma_{j}:\delta_j)=(0:1)$ for $j=r+1,\dots, s$. If $r>s$, then put $(\alpha_i:\beta _i)=(0:1)$ for $i=s+1,\dots, r$. We obtained
\[
\frac{f}{g}=\frac{a_0(\alpha_1 p+\beta_1 q)\dots(\alpha_k p+\beta_k q)}{b_0(\gamma_1 p+\delta_1 q)\dots(\gamma_k p+\delta_k q)},
\] which means that up to multiplication by a non-zero constant
 $f=(\alpha_1 p+\beta_1 q)\dots(\alpha_k p+\beta_k q)$ and
 $g=(\gamma_1 p+\delta_1 q)\dots(\gamma_k p+\delta_k q)$.

Suppose now  that $f$ and $g$ have the form as in the
conditions of this Theorem. Let us show that the rational
function $\frac{f}{g}$ is not closed. As $f=(\alpha_1 p+\beta_1 q)\dots(\alpha_k p+\beta_k q)$
and $g=(\gamma_1 p+\delta_1 q)\dots(\gamma_k p+\delta_k q)$, one has
\[
\frac{f}{g}=\frac{(\alpha_1 p+\beta_1 q)\dots(\alpha_k p+\beta_k q)}{(\gamma_1 p+\delta_1
q)\dots(\gamma_k p+\delta_k q)}=\frac{(\alpha_1 \frac{p}{q}+\beta_1 )\dots(\alpha_k \frac{p}{q}+\beta_k )}{(\gamma_1 \frac{p}{q}+\delta_1 )\dots(\gamma_k \frac{p}{q}+\delta_k )},
\]
i. e., $\frac{f}{g}=F(\frac{p}{q})$ for the rational
function
 $F(t)=\frac{(\alpha_1t+\beta_1 )\dots(\alpha_k t+\beta_k )}{(\gamma_1
 t+\delta_1 )\dots(\gamma_k t+\delta_k )}$. Since $(\alpha_i:\beta_i)\neq(\gamma_j:\delta_j)$,
 $i,j=\overline{1,k}$, we conclude that $\deg F(t)\ge 2$, which means that $\frac{f}{g}$
 is not closed (equivalently, by Theorem~\ref{th:criterion}, the pencil $\alpha f+\beta g$
 contains infinitely many reducible hypersurfaces).
\end{proof}

\begin{ex}
Let $p$ and $q$ be irreducible algebraically independent
polynomials in $\poly{n}$, $n\ge 2$. Then
$\varphi=\frac{p^l}{q^m}$ is a closed rational function for
coprime $l$ and $m$.

Indeed, suppose the converse holds. Then by
Theorem~\ref{tr:reducible} there exists irreducible polynomials
$p_1$ and $q_1$, an integer $k\ge 2$ such that
\[
\frac{p^l}{q^m}=\frac{(\alpha_1p_1+\beta_1q_1)\dots(\alpha_kp_1+\beta_kq_1)}{(\gamma_1p_1+\delta_1q_1)\dots(\gamma_kp_1+\delta_kq_1)},
\quad (\alpha_i:\beta_i)\neq (\gamma_j:\delta_j),\quad i,
j=\overline{1, k}.
\]
Since $\alpha_ip_1+\beta_iq_1$ and $\gamma_jp_1+\delta_j q_1$ are
coprime for all $i$ and $j$, it follows that
$p^l=(\alpha_1p_1+\beta_1q_1)\dots(\alpha_kp_1+\beta_kq_1)$ and
$q^m=(\gamma_1p_1+\delta_1q_1)\dots(\gamma_kp_1+\delta_kq_1)$.
Since $p$
 and $q$ are algebraically independent, as in the proof of
 Theorem~\ref{th:condition_for_closed} we conclude
  that $\alpha p_1+\beta q_1\neq \k$ for  all $(\alpha:\beta)\in \bb P^1$. Therefore,
 \(
 (\alpha_1:\beta_1)=\dots=(\alpha_k:\beta_k),
 \)
 \(
  (\gamma_1:\delta_1)=\dots=(\gamma_k:\delta_k),
 \)
  and
 \[
 p^l=a_0(\alpha_1 p_1+\beta_1 q_1)^k,\quad q^m=b_0(\gamma_1 p_1+\delta_1 q_1)^k.
 \]
 Since $p$ and $q$ are irreducible,  we obtain, up to
 multiplication by a non-zero constant,
  $\alpha_1 p_1+\beta_1 q_1=p^{l'}$ and $\gamma_1 p_1+\delta_1 q_1=q^{m'}$, i. e., $k\ge 2$
  divides both $l$ and $m$. This is impossible, since $l$ and $m$ are coprime.
  We obtained a contradiction, which proves that $\frac{p^l}{q^m}$ is a closed rational function.
 \end{ex}

\section{Products of irreducible polynomials
}

\begin{tr}
 Let $p_1,\dots, p_k \in\poly{n}$ be irreducible algebraically
independent polynomials. If $\gcd(m_{1}, m_{2}, \ldots m_{k})=1$
then the polynomial
\[
p_1^{m_{1}} p_2^{m_{2}}\dots
p_k^{m_{k}}+\lambda
\]
is irreducible for all but finitely many
$\lambda\in \k$.

\end{tr}
\begin{proof}
Show at first that the  polynomial $p_1^{m_{1}} p_2^{m_{2}}\dots
p_k^{m_{k}}$ is  closed. Suppose to the contrary  it is not closed and let
$p_1^{m_{1}} p_2^{m_{2}}\dots p_k^{m_{k}}=F(h)$ for some closed
polynomial $h$ and $F(t)\in \k[t]$, $\deg F(t)\ge 2$. Let
$F(t)=\alpha(t-\mu_1)\dots(t-\mu_s)$ be the decomposition of
$F(t)$ into linear factors. Then
\[
p_1^{m_{1}} p_2^{m_{2}}\dots
p_k^{m_{k}}=\alpha(h-\mu_1)\dots(h-\mu_s),\quad \mu\in \k,\quad
\alpha \in \k ^\star.
\]
Since the polynomials $h-\mu_i$ are closed and since we assumed
that the polynomial $p_1^{m_{1}} p_2^{m_{2}}\dots p_k^{m_{k}}$ is not closed, one
concludes that $s\ge 2$. Suppose there exists $\mu_i\neq \mu_j$,
assume without loss of generality $\mu_1\neq \mu_2$. As all $p_i$
are irreducible, it is clear that
$h-\mu_1=\alpha_1p_{i_1}^{s_{1}}\dots p_{i_m}^{s_{m}}$ and
$h-\mu_2=\alpha_2p_{j_1}^{t_{1}}\dots p_{j_r}^{t_{r}}$ for
$p_{i_1},\dots, p_{i_m},p_{j_1},\dots, p_{j_r}\in \{p_1,\dots,
p_k\}$. Since $\mu_1\neq \mu_2$, the polynomials $h-\mu_1$ and
$h-\mu_2$ are coprime. Therefore, the sets $\{p_{i_1},\dots,
p_{i_m}\}$ and $\{p_{j_1},\dots, p_{j_r}\}$ are disjoint. From
$(h-\mu_1)-(h-\mu_2)+(\mu_1-\mu_2)=0$ it follows that
\[
\alpha_1p_{i_1}^{s_{1}}\dots
p_{i_m}^{s_{m}}-\alpha_2p_{j_1}^{t_{1}}\dots
p_{j_r}^{t_{r}}+(\mu_1-\mu_2)=0,
\]
which means that the set $ \{p_1,\dots, p_k\}$ is algebraically
dependent. We obtained a contradiction. Therefore,
$\mu_1=\dots=\mu_s$ and $p_1^{m_{1}}\dots p_k^{m_{k}}=\alpha
(h-\mu_1)^s$, $s\ge 2$. From the unique factorization of the
polynomial $p_1^{m_{1}}\dots p_k^{m_{k}}$  it follows that $s|
m_{1}, \ldots , s| m_{k}$ which is impossible by our restriction
on numbers $m_{1}, \ldots , m_{k}.$ This contradiction proves that
the
 polynomial $p_1^{m_{1}}\dots p_k^{m_{k}}$ is closed. Therefore, by
  Remark~\ref{clos} the polynomial
  $p_1^{m_{1}}\dots p_k^{m_{k}} +\lambda$ is irreducible
 for all but finitely many $\lambda\in \k$.

\end{proof}
The authors are grateful to
Prof. A. Bodin who has observed on some intersection of this paper
with his preprint \cite{Bo} (in fact, the statement of Theorem 2
is equivalent to Theorem 2.2 from \cite{Bo} in zero characteristic, but the proofs of
these results are quite different).

\end{document}